\documentclass[12pt,a4paper]{article}
\usepackage[margin=2cm]{geometry}
\usepackage[english]{babel}
\usepackage[latin1]{inputenc}
\usepackage{amsmath, amsthm, amssymb}
\usepackage[pdftex]{graphicx}
\usepackage{amssymb}
\usepackage{latexsym}
\usepackage{amstext}
\usepackage{epstopdf}

\newcommand{\ad}{\mathbf{a_2}}

\newcommand{\af}{\mathbf{a_5}}
\newcommand{\afm}{\mathbf{a_{5 \frac12}}}
\newcommand{\afmH}{\mathbf{a_{5}{}_\bullet H_{20}}}
\newcommand{\au}{a_1}
\newcommand{\aud}{a_1^\dagger}
\newcommand{\aun}{\mathbf{a_1}}
\newcommand{\aut}{a_1^T}

\newcommand{\bu}{b_1}
\newcommand{\bud}{b_1^\dagger}

\newcommand{\but}{b_1^T}
\newcommand{\be}{\begin{equation}}
\newcommand{\bea}{\begin{equation}\begin{array}}
\newcommand{\beas}{\begin{equation*}\begin{array}}
\newcommand{\bef}{\begin{flalign}}
\newcommand{\befs}{\begin{flalign*}}
\newcommand{\bes}{\begin{equation*}}

\newcommand{\ct}{\mathbf{c_3}}
\newcommand{\ctm}{\mathbf{c_{3\frac12}}}
\newcommand{\ctmH}{\mathbf{c_{3}{}_\bullet H_{14}}}
\newcommand{\cu}{c_1}

\newcommand{\cut}{c_1^T}
\newcommand{\dprime}{^{\prime\prime}\mkern-1.2mu}
\newcommand{\ds}{\mathbf{d_6}}
\newcommand{\dsm}{\mathbf{d_{6 \frac12}}}
\newcommand{\dsmH}{\mathbf{d_{6}{}_\bullet H_{32}}}
\newcommand{\dsmm}{\mathbf{d_{6 \frac12 \frac12}}}
\newcommand{\dsmmH}{\mathbf{d_{6}{}_\bullet H_{32}{}_\bullet H_{44}}}
\newcommand{\ee}{\end{equation}}
\newcommand{\eea}{\end{array}\end{equation}}
\newcommand{\eeas}{\end{array}\end{equation*}}
\newcommand{\eef}{\end{flalign}}
\newcommand{\eefs}{\end{flalign*}}
\newcommand{\ees}{\end{equation*}}

\newcommand{\ekme}{\mathbf e_k^-}
\newcommand{\ekmp}{\varepsilon_k^{\mp}}

\newcommand{\ekpe}{\mathbf e_k^+}
\newcommand{\ekpm}{\varepsilon_k^{\pm}}

\newcommand{\ejpm}{\varepsilon_j^\pm}
\newcommand{\eo}{\mathbf{e_8}}

\newcommand{\es}{\mathbf{e_6}}
\newcommand{\esm}{\mathbf{e_{7 \frac12}}}
\newcommand{\esmH}{\mathbf{e_{7}{}_\bullet H_{56}}}
\newcommand{\est}{\mathbf{e_7}}

\newcommand{\fq}{\mathbf{f_4}}
\newcommand{\fqe}{\mathfrak{f}}
\newcommand{\gI}{\textbf{\large ${\mathfrak g}_{I}$}}
\newcommand{\gII}{\textbf{\large ${\mathfrak g}_{II}$}}
\newcommand{\gIII}{\textbf{\large ${\mathfrak g}_{III}$}}
\newcommand{\gIIIm}{\textbf{\large ${\mathfrak g}_{III \frac{1}{2}}$}}
\newcommand{\gIV}{\textbf{\large ${\mathfrak g}_{IV}$}}
\newcommand{\gd}{\mathbf{g_2}}

\newcommand{\jdot}{\!\cdot\!}

\newcommand{\jotn}{\mathbf{J_3^n}}
\newcommand{\jobtn}{\mathbf{\overline J_3^{\raisebox{-2 pt}{\scriptsize \textbf n}}}}
\newcommand{\jotu}{\mathbf{J_3^1}}

\newcommand{\jotd}{\mathbf{J_3^2}}

\newcommand{\joto}{\mathbf{J_3^8}}
\newcommand{\jobto}{\mathbf{\overline J_3^{\raisebox{-2 pt}{\scriptsize \textbf 8}}}}
\newcommand{\jots}{\mathbf{J_3^6}}
\newcommand{\jobts}{\mathbf{\overline J_3^6}}

\newcommand{\Lvsp}{L_{\mathbf s^+}}
\newcommand{\Lvsm}{L_{\mathbf s^-}}
\newcommand{\Lvspm}{L_{\mathbf s^\pm}}

\newcommand{\Lvrp}{L_{\mathbf r^+}}
\newcommand{\Lvrm}{L_{\mathbf r^-}}
\newcommand{\Lvrpm}{L_{\mathbf r^\pm}}

\newcommand{\mep}{\star}

\newcommand{\nbf}{\mathbf{n}}
\newcommand{\ncirc}{\underline{\circ}\ }                 
\newcommand{\oo}{\textbf{\large $\mathfrak C$}}

\newcommand{\rep}{\mathbf \varrho}

\newcommand{\rmp}{\rho^\mp}

\newcommand{\rpm}{\rho^\pm}
\newcommand{\rr}{\mathbf{R}}
\newcommand{\sex}{\mathbb S} 

\newcommand{\vrm}{\mathbf r^-}
\newcommand{\vrp}{\mathbf r^+}
\newcommand{\vrpm}{\mathbf r^\pm}
\newcommand{\vsm}{\mathbf s^-}
\newcommand{\vsp}{\mathbf s^+}
\newcommand{\vspm}{\mathbf s^\pm}

\newcommand{\vx}{\mathbf x}

\newcommand{\vy}{\mathbf y}

\newcommand{\xs}{x^\#}

\numberwithin{equation}{section}
\begin{document}
%
\begin{titlepage}
\begin{center}

\hfill DFPD/2015/TH/10

\vskip 2cm


{\huge{\bf Sextonions, Zorn Matrices, and $\esm$}}

\vskip 2.5cm

{\large{\bf Alessio Marrani\,$^{1,2}$ and  Piero Truini\,$^3$}}

\vskip 20pt

{\it $^1$ Centro Studi e Ricerche ``Enrico Fermi'',\\
Via Panisperna 89A, I-00184, Roma, Italy \vskip 5pt }

\vskip 10pt

{\it $^2$ Dipartimento di Fisica e Astronomia ``Galileo Galilei'', \\Universit\`a di Padova,\\and INFN, sezione di Padova,\\Via Marzolo 8, I-35131 Padova, Italy}\\\vskip 5pt

\texttt{Alessio.Marrani@pd.infn.it}\\

    \vspace{10pt}

{\it ${}^3$ Dipartimento di Fisica, Universit\` a degli Studi\\
and INFN, sezione di Genova,\\
via Dodecaneso 33, I-16146 Genova,  Italy}\\\vskip 5pt
\texttt{truini@ge.infn.it}

    \vspace{10pt}

\end{center}

\vskip 2.2cm

\begin{center} {\bf ABSTRACT}\\[3ex]\end{center}
By exploiting suitably constrained Zorn matrices, we present a new
construction of the algebra of sextonions (over the algebraically closed
field $\mathbb{C}$). This allows for an explicit construction, in terms of
Jordan pairs, of the non-semisimple Lie algebra $\esm$,
intermediate between $\est$ and $\eo$, as well as of
all Lie algebras occurring in the sextonionic row and column of the extended
Freudenthal Magic Square.

\vskip 0.2cm

MSC(2010) numbers: 17B10, 17B25, 17B45.

\vskip 0.1cm

Keywords : exceptional Lie algebras, intermediate algebras, sextonions, Zorn matrices.

\end{titlepage}

\newpage \setcounter{page}{1} \numberwithin{equation}{section}


\section{Introduction}

The field of composition, non-associative algebras, and related Lie
algebras, underwent a series of interesting developments in recent times.

In \cite{D-1} Deligne proposed dimension formulas for the exceptional series
of complex simple Lie algebras, whose parametrization in terms of the dual
Coxeter number was exploited further in \cite{CD-1} by Cohen and de Man (see
also \cite{D-2}). Landsberg and Manivel subsequently pointed out the
relation between the dimension formulas and the dimensions of the
composition algebras themselves in \cite{LM-1}. In \cite{D-1,CD-1} it was
observed that all parameter values determining integer outputs in the
dimension formulas were already accounted for by the known normed division
algebras, with essentially one exception, intriguingly corresponding to a
would-be composition algebra of dimension six, sitting between the
quaternions and octonions.

This algebra, whose elements were named \textit{sextonions}, was recently
studied by Westbury in \cite{W-1}, who pointed out the related existence of
a whole new row in the Freudenthal Magic Square. Actually, the
six-dimensional algebra of sextonions had been observed earlier as a
curiosity; indeed, it was explicitly constructed in \cite{Kle}. Moreover, it
was used in \cite{Jeu} to study the conjugacy classes in the smallest
exceptional Lie algebra $\mathbf{g}_{2}$ in characteristics other than $2$
or $3$. The sextonions were also constructed in \cite{Rac} (\textit{cfr.}
Th. 5 therein), and proved to be a maximal subalgebra of the split octonions.

In \cite{LM-2}, Landsberg and Manivel \textquotedblleft filled in the hole"
in the exceptional series of Lie algebras, observed by Cvitanovic, Deligne,
Cohen and de Man, showing that sextonions, through the \textit{triality}
construction of \cite{LM-1}, give rise to a non-simple \textit{intermediate}
exceptional Lie algebra, named $\esm$, between $\est$ and $\eo$, satisfying
some of the decomposition and dimension formulas of the exceptional simple
Lie algebras \cite{D-1,CD-1,D-2,LM-1,LM-3}.

More recently, such a $190$-dimensional Lie algebra $\esm$ was also found by
Mkrtchyan in the study of the Vogel plane \cite{MK-1}, in the context of the
analysis of the \textit{universal} Vogel Lie algebra \cite{V-1}.\medskip

By the Hurwitz Theorem \cite{H-1}, the real normed division algebras are the
real numbers $\mathbb{R}$, the complex numbers $\mathbb{C}$, the quaternions
$\mathbb{H}$ and the octonions $\mathfrak{C}$ (Cayley numbers). Each algebra
can be constructed from the previous one by the so-called \textit{%
Cayley-Dickson doubling procedure }\cite{Di-1,Sch-1}.

All these algebras can be complexified to give complex algebras. These
complex algebras respectively are

$\mathbb{R}\otimes \mathbb{C}$ $=$ $\mathbb{C}$, $\mathbb{C}\otimes \mathbb{C%
}$ $=$ $\mathbb{C}\oplus \mathbb{C}$, $\mathbb{H}\otimes \mathbb{C}$ $=$ $%
M_{2}(\mathbb{C})$, $\mathfrak{C}\otimes \mathbb{C}$ ($M_{2}$ denoting a $%
2\times 2$ matrix). The three complex algebras other than $\mathbb{C}$ have
a second real form, denoted $\mathbb{C}_{s}$, $\mathbb{H}_{s}$ and $%
\mathfrak{C}_{s}$, with the following isomorphisms holding : $\mathbb{C}_{s}=%
\mathbb{R}\oplus \mathbb{R}$ and $\mathbb{H}_{s}=$ $M_{2}(\mathbb{R})$. The
normed division algebras are called the \textit{compact} forms and the
aforementioned second real form is called the \textit{split} real form. It
is worth pointing out that split real forms are composition algebras but
they are not division algebras.

On the field $\mathbb{R}$, the sextonions only exist in split form $\mathbb{S%
}_{s}$, and they are intermediate between

the split quaternions $\mathbb{H}_{s}$ and the split octonions $\mathfrak{C}%
_{s}$ :%
\begin{equation}
\mathbb{H}_{s}\subset \mathbb{S}_{s}\subset \mathfrak{C}_{s}.
\end{equation}%
Note that $\mathbb{S}_{s}$ does not contain the divisional quaternions $%
\mathbb{H}$; see App. A.\bigskip

Nowadays, exceptional Lie algebras have a long-standing history of
applications to physics (see \textit{e.g.} \cite{1-E}, \cite{Toppan}, \cite{1}%
-\nocite{2,3,4,5,6,7}\cite{8} for a partial list of results and Refs.). The
relevance of compact exceptional Lie algebras (and groups) in realizing
grand unification gauge theories and consistent string theories is well
recognized. Similarly, the relevance of non-compact real forms for the
construction of locally supersymmetric theories of gravity is well
appreciated. Other frameworks include sigma models based on quotients of
exceptional Lie groups, which are of interest for string theory and
conformal field theories, as well. It is here worth pointing out that that
the analysis of quantum criticality in Ising chains and the structure of
magnetic materials such as Cobalt Niobate has also recently (and strikingly)
turned out to be related to exceptional Lie algebras of type $E$ (see
\textit{e.g.} \cite{2-E} and \cite{3-E}, respectively). Moreover,
exceptional Lie algebras occur in models of confinement in non-Abelian gauge
theories (for instance, \textit{cfr.} \cite{4}), as well as in a striking
relation between cryptography and black hole physics, recently discovered
\cite{5-E-1}-\nocite{5-E-2,5-E-3,6-E,7-E}\cite{8-E}. It should also be recalled that
fascinating exceptional algebraic structures arise in the description of the
Attractor Mechanism for black holes in Maxwell-Einstein supergravity
theories \cite{9-E-1}-\nocite{9-E-2,9-E-3,9-E-4,9-E-5,FG}\cite{AM}, such as the so-called magic exceptional
supergravity \cite{10-E-1}-\nocite{10-E-2,10-E-3,10-E-4}\cite{11-E}.

In this context, the aforementioned, intermediate $190$-dimensional Lie
algebra $\esm$ is quite novel, and applications to physics are still under
investigation, even though recent studies (\textit{cfr. e.g.} \cite{Devchand}%
) intriguingly seem to connect sextonions to theories \textit{beyond}
eleven-dimensional $M$-theory. It should also be mentioned that $\esm$ can
be regarded as a \textit{Freudenthal triple system} over the exceptional
Albert algebra $J_{3}^{\mathbb{O}}$, along with its automorphism algebra\ $%
\est$ and an extra $\est$-singlet generator, acting on the Freudenthal
triple system as the multiplication by a scalar; for details on the
applications of Freudenthal triple systems to the study of black hole
attractors in four dimensions, \textit{cfr. e.g.} \cite{FG, BFGM1,
Small-Orbits, Small-Orbits-maths}, and Refs. therein.\bigskip

In the present paper we will apply the formal machinery introduced in \cite%
{T-1} and \cite{T-2}, as well as an explicit realization of the sextonions
(over the algebraically closed field $\mathbb{C}$), in order to explicitly
construct the non-semisimple Lie algebra $\esm$, as well as all algebras
occurring in the sextonionic row of the \textit{extended} Freudenthal Magic
Square \cite{W-1,LM-2}, in terms of Jordan pairs.

\bigskip

The plan of the paper is as follows.

In Sec. \ref{sec:octonions} we provide a realization of the sextonions in
terms of nilpotents constructed from the traceless octonions, and recall their representation in terms of suitably constrained Zorn matrices.

The intermediate exceptional algebra $\esm$ is then considered in Sec. \ref%
{jazz}, which focuses on the construction (then developed in Secs. \ref%
{Blues1} and \ref{Blues2}) of the sextonionic row and column of the extended
Magic Square, by exploiting Jordan pairs for the sextonionic rank-3 Jordan
algebra.

The action of $\mathbf{g}_{2}=Der(\mathfrak{C})$ on the Zorn matrices is
recalled in Sec. \ref{g2d}, and exploited in Sec. \ref{sec:d(s)} to
determine the derivations of the sextonions, $Der(\sex)$.

Then, in Secs. \ref{Blues1} and \ref{Blues2} the explicit construction of
the intermediate algebras $\mathbf{c}_{3\frac{1}{2}}$ (which analogously
holds for $\mathbf{a}_{5\frac{1}{2}}$ and $\mathbf{d}_{6\frac{1}{2}}$) and $%
\esm$ is presented.

The paper is concluded by App. A, in which we prove that, on the field $%
\mathbb{R}$, $\mathbb{S}_{s}$ does not contain the divisional quaternions $%
\mathbb{H}$.

\section{Sextonions and their Nilpotent Realization}\label{sec:octonions}
\begin{figure}
\begin{center}
\includegraphics{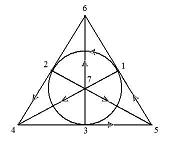}
\caption{Fano diagram for the octonions' products}\label{fig:fano}
\end{center}
\end{figure}

The algebra of sextonions is a six dimensional subalgebra $\sex$ of the octonions. As mentioned above, we denote by
$\oo$ the algebra of the octonions over the complex field  $\mathbb{C}$, whose multiplication rule goes according to the Fano diagram in Figure \ref{fig:fano}.

If $a \in \oo$ we write $a = a_0 + \sum_{j=1}^7{a_j u_j}$ where $a_j \in \mathbb{C}$ for $j = 1, \dots , 7$ and $u_j$ for $j = 1, \dots , 7$ denote the octonion imaginary units. We denote by $i$ the the imaginary unit in $\mathbb{C}$.

We introduce 2 idempotent elements:
$$\rpm = \frac{1}{2}(1 \pm i u_7) $$
and 6 nilpotent elements :
$$\ekpm = \rpm u_k \quad , \qquad k = 1,2,3$$
One can readily check that \cite{T-1}:
\begin{equation}
\begin{array}{ll}
 & (\rpm)^2 = \rpm \quad , \quad \rpm \rmp = 0 \\ \\
 & \rpm \ekpm = \ekpm \rmp = \ekpm \\ \\
 & \rmp \ekpm = \ekpm \rpm = 0 \\ \\
 & (\ekpm)^2 = 0  \\ \\
& \ekpm \varepsilon_{k+1}^\pm = - \varepsilon_{k+1}^\pm \ekpm  = \varepsilon_{k+2}^\mp \qquad \text{(indices modulo 3)} \\ \\
& \ejpm \ekmp = 0 \qquad j \ne k \\ \\
& \ekpm \ekmp = - \rpm
\end{array}\end{equation}

We can write $a \in \oo$ as $a = \alpha_0^+ \rho^+ +\alpha_0^- \rho^- + \alpha_k^+
\varepsilon_k^+ +\alpha_k^- \varepsilon_k^-$.

The subalgebra $\sex \in \oo$ generated by $\rpm , \varepsilon_1^\pm , \varepsilon_2^+ , \varepsilon_3^-$ (namely $a \in \sex$ \textit{iff} $\alpha_2^- = \alpha_3^+ = 0$)  provides an explicit realization of the sectonions. The existence of
the non-divisional sextonionic elements can be easily understood. Indeed, in order to
construct divisional sextonions, one would need to combine a nilpotent with
its complex conjugate; but, as given by the above construction, this is not
possible for $\varepsilon _{2}^{+}$ nor for $\varepsilon _{3}^{-}$.

Octonions can be represented by Zorn matrices \cite{zorn}. After the treatment of Sec. 3 of \cite{T-2}, we can represent the sextonions as a Zorn matrix, as long as  $A^+$ and $A^-$ are  $\mathbb{C}^3$-vectors of the type
$$A^+ = (a^+,c^+,0) \qquad \text{and} \qquad A^- = (a^-,0,c^-)$$
Notice that $A^+$ and $A^-$ lie on orthogonal $\mathbb{C}^3$-planes sharing the line along the first component.

\section{$\esm$}\label{jazz}

In recent papers \cite{T-1, T-2} a unifying view of all exceptional Lie algebras in terms of $\ad$ subalgebras and {\it Jordan Pairs} has been presented, and a {\it Zorn-matrix-like} representation of these algebras has been introduced.

\begin{figure}
\begin{center}
\includegraphics{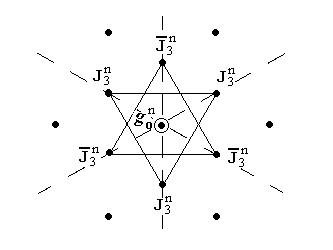}
\caption{A unifying view of the roots of exceptional Lie algebras through \textit{Jordan pairs} \cite{T-1}. For $\mathbf{n}=8$, the root diagram of $\mathbf{e}_{8}$ is obtained.}\label{fig:rootdiagram}
\end{center}
\end{figure}

The root diagram related to this view is shown in Figure \ref{fig:rootdiagram}, where the roots of the
exceptional Lie algebras are projected on a complex $\mathbf{su(3)}=\mathbf{a}_{2}$ plane,
recognizable by the dots forming the external hexagon, and it exhibits the
\textit{Jordan pair} content of each exceptional Lie algebra. There are
three Jordan pairs $(\jotn,\jobtn)$, each of which lies on an axis
symmetrically with respect to the center of the diagram. Each pair doubles
a simple Jordan algebra of rank $3$, $\jotn$, with involution - the
conjugate representation $\jobtn$, which is the algebra of $3\times 3$
Hermitian matrices over $\mathbb{A}$, where $\mathbb{A}=\mathbb{R},\mathbb{C},\mathbb{H},\mathfrak{C}$ for $\mathbf{n}=$dim$_{\mathbb{R}}\mathbb{A}=1,2,4,8$ respectively, stands for real, complex, quaternion, octonion algebras, the four normed division algebras according to Hurwitz's Theorem; see
\textit{e.g.} \cite{McCrimmon}.
Exceptional Lie algebras $\fq$, $\es$, $\est$, $\eo$ are obtained for $\nbf%
=1,2,4,8$, respectively. $\gd$ (corresponding to $\nbf
=-2/3$) can be also represented in the same way, with
the Jordan algebra reduced to a single element. For further detail, \textit{cfr.} \cite{T-1}.

\begin{figure}
\begin{center}
\includegraphics{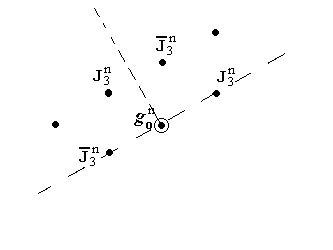}
\caption{Root diagram of $\esm$ (for $\mathbf{n}=8$)}\label{fig:rootesm}
\end{center}
\end{figure}

We expand that view in this paper to include $\esm$ \cite{LM-2}, a Lie subalgebra of $\mathbf{e}_8$ of dimension 190. If we consider the $\mathbf{e}_8$ root diagram (obtained in Figure \ref{fig:rootdiagram} for $\mathbf{n}=8$),
then the sub-diagram of $\esm$ is shown in Figure \ref{fig:rootesm}, (for $\mathbf{n}=8$, as well).

%
%

In general, one can do the same for all algebras in the fourth and third row of the Magic square, \cite{tits2} \cite{freu1}, that we denote by $\gIV$ and $\gIII$ respectively (see table \ref{ms34}). In this way, the algebras in the intermediate (fourth) row of the extended Magic Square \cite{W-1, LM-2} are explicitly constructed in terms of Jordan pairs.
\begin{table}[h!]
\begin{center}
\begin{tabular}{|c||c|c|c|c|c|c|c|}
\hline
$\mathbf{n}$ & $1$ & $2$ & $4$ & $8$ \\ \hline
\rule[-1mm]{0mm}{6mm} $\gIII$ & $\ct$ & $\af$ & $\ds$ & $\est$ \\ \hline
\rule[-1mm]{0mm}{6mm} $\gIV$ & $\fq$ & $\es$ & $\est$ & $\eo$ \\ \hline
\end{tabular}
\end{center}

\caption{Third and fourth row of the magic square \label{ms34}}
\end{table}

\begin{figure}
\begin{center}
\includegraphics{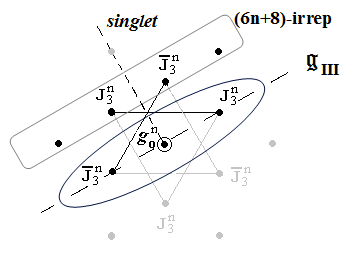}
\caption{ Diagram of $\gIIIm$}\label{fig:gIIIm}
\end{center}
\end{figure}

We get a subalgebra of $\gIV$, that we denote here by $\gIIIm$, given by $\gIII$ plus a $(6n+8)$-dimensional irreducible representation of $\gIII$  plus a $\gIII$-singlet, as shown in Figure \ref{fig:gIIIm} \footnote{There are some variations on the definition of {\it intermediate} algebra, \cite{W-1}
\cite{LM-2}, based on the grading induced by an highest root. Our realisation of $Der(\sex)$ and $\esm$ corresponds to the algebra denoted by $\bf{{\mathfrak g}\dprime}$ in the Introduction of \cite{LM-2}.}.

In particular, the irreps. of $\gIII$ are symplectic (\textit{i.e.}, they admit a skew-symmetric invariant form), and they have complex dimension $6n+8=14, 20, 32, 56$ for $n = 1, 2, 4, 8$ respectively; the algebras  $\gIIIm$ are their corresponding Heisenberg algebras (denoted by $\mathbf{H}$) through such an invariant tensor, \cite{LM-2}, $\ctm = \ctmH$, $\afm = \afmH$, $\dsm = \dsmH$, $\esm = \esmH$, of complex dimension $36, 56, 99, 190$.\\ 

Let us here present a brief account of the \textit{Jordan pairs} for
sextonions $\mathfrak{\sex}$ by means of suitable embeddings. We start with the
maximal, non-symmetric embedding:%
\begin{eqnarray}
\mathbf{e}_{7} &\supset &\mathbf{a}_{2}\oplus \mathbf{a}_{5} \\
\mathbf{133} &=&\left( \mathbf{8},\mathbf{1}\right) +\left( \mathbf{1},%
\mathbf{35}\right) +\left( \mathbf{3},\mathbf{15}\right) +\left( \overline{%
\mathbf{3}},\overline{\mathbf{15}}\right)  \\
\mathbf{56} &=&\left( \mathbf{3},\mathbf{6}\right) +\left( \overline{\mathbf{%
3}},\overline{\mathbf{6}}\right) +\left( \mathbf{1},\mathbf{20}\right) ,
\end{eqnarray}%
implying that:%
\begin{equation}
\mathbf{e}_{7}\ltimes \mathbf{56}\supset \left[ \mathbf{a}_{2}\oplus \left(
\mathbf{a}_{5}\ltimes \mathbf{20}\right) \right] \ltimes \left( \mathbf{3},%
\mathbf{15}+\mathbf{6}\right) +\left( \overline{\mathbf{3}},\overline{%
\mathbf{15}}+\overline{\mathbf{6}}\right) .\label{decc}
\end{equation}%
Thus, the \textit{Jordan pairs} for the sextonionic Jordan algebra of rank
3, $\mathbf{J}_{3}^{\mathbf{n}=6}$, are given by $\left( \mathbf{3},\mathbf{%
15}+\mathbf{6}\right) +\left( \overline{\mathbf{3}},\overline{\mathbf{15}}+%
\overline{\mathbf{6}}\right) $ in (\ref{decc}).

In order to reconstruct the extended Magic Square \cite{W-1, LM-2}, one needs also to add the extra column shown in table \ref{ms6}, \cite{LM-2},  
where a further algebra $\dsmm = \dsmmH$   
is introduced.

\begin{table}[h!]
\begin{center}
\begin{tabular}{|c||c|}
\hline $\mathbf{n}$ & $6$  \\ \hline
\rule[-1mm]{0mm}{6mm} $\gI$ & $\ctm$ \\ \hline 
\rule[-1mm]{0mm}{6mm} $\gII$ & $\afm$ \\ \hline 
\rule[-1mm]{0mm}{6mm} $\gIII$ & $\dsm$ \\ \hline
\rule[-1mm]{0mm}{6mm} $\gIIIm$ & $\dsmm$ \\ \hline 
\rule[-1mm]{0mm}{6mm} $\gIV$ & $\esm$ \\ \hline
\end{tabular}
\end{center}

\caption{Sixth column of the magic square \label{ms6}}
\end{table}

This column corresponds to the Jordan algebra that we denote by
$\jots$  
of $3 \times 3$ Hermitian matrices over the sextonions. The \emph{new} element $\dsmmH$ can be easily seen in the diagram of figure \ref{fig:gIIIm} for $n=6$: $\mathbf{g_0^6} = \afm$ is the reduced structure algebra of $\jots$,\,  $\gIII = \dsm$ the super-structure algebra of $\jots$ and finally $\dsmm = \dsm \mathbf{{}_\bullet H_{44}} = \dsmmH$. Notice that the $44$-dimensional representation of $\dsm$ is made of $\jots \oplus \jobts \oplus 2$.  
Finally, the algebra $\esm$ at the end of the column is viewed as in the diagram of Figure \ref{fig:rootdiagram} for $n = 6$, with $\mathbf{g_0^6} = \afm$ and the subalgebra $\est$ represented by the same diagram for $n = 4$.

This completes the explicit construction of the relevant rows and columns (pertaining to the sextonions) of the extended Magic Square\footnote{It is once again worth stressing that in the present investigation, as well as in the previous papers \cite{T-1, T-2}, we only consider complex forms of the Lie algebras.} \cite{W-1, LM-2}.
%
%

\section{$\gd$ action on Zorn matrices}\label{g2d}

In our previous paper \cite{T-2}, we have introduced the following adjoint representation $\rep$ of the Lie algebra $\gd$:
\begin{equation}
\left[
\begin{array}{cc}
a &  A^+ \\
A^- & 0
\end{array} \right]
\label{gdm2} \end{equation}
where $a \in \ad$, $A^+,\  A^- \in \mathbb{C}^3$, viewed as column and row vector respectively.\\
The commutator of two such matrices reads \cite{T-2}:
\begin{equation}
\begin{array}{c}
\left [ \left[ \begin{array}{cc} a & A^+ \\ A^- & 0
\end{array}\right]
 ,
\left [ \begin{array}{cc} b & B^+ \\ B^- & 0
\end{array}\right]\right] \\ =
\left [\begin{array}{cc} [a, b] + A^+ \circ  B^- - B^+\circ A^- &
a B^+  - b A^+ + 2 A^- \wedge B^- \\
A^- b - B^-a + 2 A^+ \wedge B^+ & 0 \end{array} \right]
\label{comm}
\end{array}
\end{equation}
where
\be
A^+ \circ B^- = t(A^+ B^- )  I - t(I) A^+ B^-
\label{circ}
\ee
(with standard matrix products of row and column vectors and with $I$ denoting the $3\times 3$ identity matrix); $A
\wedge B$ is the standard vector product of $A$ and $B$, and $t(a)$ denotes the trace of $a$.\\

\begin{figure}
\begin{center}
\includegraphics{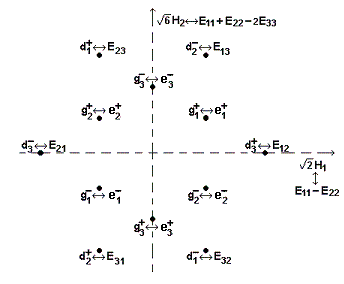}
\caption{Diagram of $\gd$ with corresponding generators and {\it matrix-like} elements}\label{fig:g2m}
\end{center}
\end{figure}

The $\gd$ generators are \cite{T-1}:
\be
\begin{array}{l}
\rep (d^\pm_k) = E_{k\pm 1 \  k\pm 2} \quad \text{(mod 3)} \ , \ k = 1,2,3\\
\rep(\sqrt{2} H_1) = E_{11} - E_{22} \qquad \rep(\sqrt{6} H_2) = E_{11} + E_{22} - 2 E_{33} \\
\rep(g^+_k) = E_{k 4} := \ekpe \quad \rep(g^-_k) = E_{4 k} := \ekme \quad  , \ k = 1,2,3
\end{array}
\ee
where $E_{ij}$ denotes the matrix with all zero elements except a $1$ in the $\{ij\}$ position: $(E_{ij})_{k\ell} = \delta_{ik}\delta_{j\ell}$ and $\ekpe$ are the standard basis vectors of $\mathbb{C}^3$ ($\ekme$ are their transpose). The correspondence with the roots of $\gd$ is shown in Figure \ref{fig:g2m}.

\section{Derivations of $\sex$ \label{sec:d(s)}}

We now use the representation $\rep$ to get a representation of the Lie algebra of $Der(\sex)$, which indeed is a non-reductive subalgebra of $\mathbf{g}_{2}=Der(\mathfrak{C})$.

It was shown in \cite{W-1} that the map from the subalgebra of derivations of $\oo$ preserving $\sex$, that we here denote by $Der_\oo (\sex)$, to $Der(\sex)$ is surjective with one-dimensional kernel; the corresponding statement at the level of automorphism group was made in \cite{LM-2}.

Within our formalism, this result is achieved by restricting $\rep(\gd)$ to the matrices that preserve $\sex$. One easily gets:
\begin{equation}
\left[
\begin{array}{cc}
a &  S^+ \\
S^- & 0
\end{array} \right]  :  a =
\left(\begin{array}{ccc}
a_{11} &  0 & a_{13} \\
a_{21} &  a_{22} & a_{23} \\
0 &  0 & a_{33} \\
\end{array} \right) \ , \ S^+ =
\left(\begin{array}{ccc}
s^+_1  \\
s^+_2 \\
0 \\
\end{array} \right)  \ , \  S^- = (s^-_1,0,s^-_3)
\label{sm} \end{equation}

We also realize very easily that the generator corresponding to $d_1^+$, namely the element $E_{23}$ in $\rep(\gd)$, acts trivially on $\sex$, hence it can be set to $0$. The commutator (\ref{comm}) must be modified accordingly, by setting the $\{ 23 \}$ element of $a$ equivalent to zero, that is by replacing the standard matrix product of two matrices
\be
 a =
\left(\begin{array}{ccc}
a_{11} &  0 & a_{13} \\
a_{21} &  a_{22} & 0 \\
0 &  0 & a_{33} \\
\end{array} \right) \ , \   b =
\left(\begin{array}{ccc}
b_{11} &  0 & b_{13} \\
b_{21} &  b_{22} & 0 \\
0 &  0 & b_{33} \\
\end{array} \right)
\ee
with the new product
\be
a \centerdot b = ab - E_{22}\ ab\ E_{33}
\label{nmpu}
\ee
and the product $S^+ \circ S^-$ with
\be
S^+ \ncirc S^- = t(S^+ S^- )  I - t(I) (S^+ S^- - E_{22} S^+ S^- E_{33})
\label{nmpd}
\ee
We thus have $Der(\sex) = \aun \oplus \mathbb{C} \oplus V_4$, where $V_4$ is a 4-dimensional\footnote{%
This representation also characterizes $\mathbf{a}_{1}$ as the smallest Lie
group \textquotedblleft of type $E_{7}$" \cite{Brown}, and it pertains to
the so-called $T^{3}$ model of $N=2$, $D=4$ supergravity.} (spin-$3/2$) irreducible representation of $\aun$ (as confirmed by the entry in the first column, fourth row in the extended Magic Square; \textit{cfr. e.g.} \cite{LM-2}). The corresponding root diagram is shown in Figure \ref{fig:derS}, where we have also included the axes corresponding  the linear span of the Cartan generators, represented by the matrices:

\begin{figure}
\begin{center}
\includegraphics{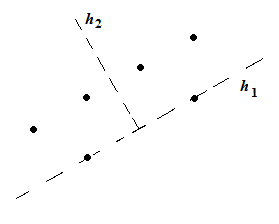}
\caption{Root diagram of $Der(\sex)$}\label{fig:derS}
\end{center}
\end{figure}

\begin{equation}
\left[
\begin{array}{cc}
h_{1,2} &  0 \\
0 & 0
\end{array} \right]  :  h_1 =
\left(\begin{array}{ccc}
-2 &  0 & 0 \\
0 &  1 & 0 \\
0 &  0 & 1 \\
\end{array} \right) \ , \ h_2 =
\left(\begin{array}{ccc}
0 &  0 & 0 \\
0 &  1 & 0 \\
0 &  0 & -1 \\
\end{array} \right)
\label{cart} \end{equation}

\noindent {\bf Proposition \ref{sec:d(s)}.1} : The algebra spanned by the generators corresponding to the roots in Fig.6 is a Lie algebra. \\
{\bf Proof :}
By looking at the diagram in Fig.5 these generators are $d^-_2, d^-_3, g^+_2, g^-_3$ spanning a subspace $L_1$ of $\gd$, plus the generators $g^\pm_1, h_1, h_2$ spanning the Lie subalgebra $L_0 := a_1\oplus \mathbb{C}$.
We have $[L_0,L_0] \subset L_0$, $[L_0,L_1] \subset L_1$, $[L_1,L_1] \subset L_2 \sim 0$, where $L_2$ is the span of $d^+_1$. The notation is that of the grading with respect to $h_2$.

 We consider the $\gd$ commutation relations among these generators and identify $d^+_1 \sim 0$. We only need to prove that the Jacobi identity is consistent with this identification.
Let $X,Y,Z \in L_0 \oplus L_1$, then consistency must be checked in only two cases (up to cyclic permutation):

\textbf{case 1:} $[X,Y] \propto d^+_1$;

\textbf{case 2:} $[[X,Y],Z] \propto d^+_1$.

\textbf{Case 1:}  Consistency requires $[[Y,Z],X] + [[Z,X],Y] \sim 0$. This is true if $[d^+_1, Z] = 0$, since it is true in $\gd$. On the other hand, if $[d^+_1, Z] \neq 0$ then $Z \propto h_2$ and $[Z,X] = \lambda X$ , $[Z,Y] = \lambda Y$ ,since $X,Y$ must be in $L_1$ by hypothesis. Therefore $[[Y,Z],X] + [[Z,X],Y] = 2 \lambda [X,Y] \propto d^+_1 \sim 0$.

\textbf{Case 2:} Both $[X,Y]$ and $Z$ must be in $L_1$. In particular either $X$ or $Y$ must be in $L_1$. Suppose $X \in L_1$. Then $[Y,Z]\in L_1$ hence we have both $[X,Z] \sim 0$ and $[[Y,Z],X] \sim 0$. Similarly if $Y \in L_1$.
\\

This concludes the proof $\blacksquare $

\section{$\mathbf{n}=1$ : Matrix representation of $\ctm$}\label{Blues1}

We denote by a dot the Jordan product $x\jdot y = \frac12(xy+yx)$ and by $ t()$ the ordinary trace of $3\times 3$ matrices. We also set $t(x,y) := t(x\jdot y)$. For $\jotu$ and $\jotd$,  obviously $t(x,y) = t(x y)$.

We use in this section the representation $\rep$ of $\fq$ in the form of a matrix introduced in \cite{T-2}, restricted to the subalgebra $\ctm$:

\begin{equation}
\rep(\fqe) = \left ( \begin{array}{cc} a\otimes I + I \otimes \au & \vsp \\ \vsm & -I\otimes \aut
\end{array}\right)
\label{mfq}
\end{equation}
where
\begin{equation}
a =
\left(\begin{array}{ccc}
a_{11} &  0 & a_{13} \\
a_{21} &  a_{22} & 0 \\
0 &  0 & a_{33} \\
\end{array} \right) \ , \ t(a)=0\ , \ \vsp =
\left(\begin{array}{ccc}
s^+_1  \\
s^+_2 \\
0 \\
\end{array} \right)  \ , \  \vsm = (s^-_1,0,s^-_3)
\label{sm2}
\ee
and $\au \in \ad$, $\aut$ is the transpose of $\au$,  $I$ is the $3\times 3$ identity matrix, $s_i^\pm \in \jotu  \ , \quad i=1,2,3$.

The commutator is set to be:

\be
\begin{array}{c}
\left[
\left ( \begin{array}{cc} a\otimes I + I \otimes \au & \vsp \\ \vsm & -I\otimes \aut
\end{array}\right) ,
\left ( \begin{array}{cc} b\otimes I + I \otimes \bu &\vrp \\ \vrm & -I\otimes \but
\end{array}\right) \right] \\  \\ :=
\left (\begin{array}{cc} C_{11} & C_{12}\\
C_{21} & C_{22}
 \end{array} \right) \hfill
\label{fqcom}
\end{array}
\ee

where, denoting by $[a\centerdot b]$ the commutator with respect to the product (\ref{nmpu})
\be [a\centerdot b] = a\centerdot b - b\centerdot a = [a,b] - E_{22}[a,b] E_{33} \ee, it holds that:

\be
\begin{array}{ll}
C_{11} &= [a\centerdot b] \otimes I + I \otimes [\au,\bu] + \vsp \diamond \vrm - \vrp \diamond \vsm \\ \\
C_{12} &=  (a \otimes I) \vrp -  (b \otimes I) \vsp + (I \otimes \au) \vrp + \vrp (I \otimes \aut) \\
 &\phantom{:=} - (I \otimes \bu) \vsp - \vsp (I \otimes \but) +  \vsm \times \vrm \\ \\
C_{21} &= - \vrm (a \otimes I)  +  \vsm (b \otimes I) - (I \otimes \aut) \vrm - \vrm (I \otimes \au)  \\
&\phantom{:=} + (I \otimes \but) \vsm + \vsm (I \otimes \bu) +  \vsp \times \vrp \\ \\
C_{22} &=  I \otimes [\aut,\but] + \vsm \bullet \vrp - \vrm \bullet \vsp
\end{array}
 \label{comrel}
\ee
with the following definitions (summing over repeated indices) :
\be
\begin{array}{ll}
\vsp \diamond \vrm &:= \left(\frac13 t(s^+_1, r^-_1) I - (1-(E_{23})_{ij}) t(s^+_i,r^-_j) E_{ij} \right) \otimes I +\\
&\phantom{:=} I \otimes \left(\frac13 t(s^+_1, r^-_1) I - s^+_1 r^-_1 \right) \\ \\
\vsm \bullet \vrp &:= I \otimes (\frac13 t(s^-_1,r^+_1) I -  s^-_1 r^+_1) \\  \\
(\vspm \times \vrpm)_i &:= \epsilon_{ijk}[s_j^\pm r_k^\pm + r_k^\pm s_j^\pm -s_j^\pm t(r_k^\pm) - r_k^\pm t(s_j^\pm) \\
&\phantom{:=}- (t(s_j^\pm, r_k^\pm) - t(s_j^\pm) t( r_k^\pm)) I]  \\
&:= \epsilon_{ijk} (s_j^\pm \# r_k^\pm)
\end{array}
 \label{not1}
\ee
Notice that:
\begin{enumerate}
\item $s \in \jotu$ is a symmetric complex matrix;
\item writing $\vsp \diamond \vrm := c \otimes I + I\otimes \cu$ we have that both $c$ and $\cu$ are traceless hence $c$ ia a matrix like $a$ in (\ref{sm}), $\cu \in \ad$ and $\vrm \bullet \vsp = I\otimes \cut$
\item terms like $(I \otimes \au) \vrp + \vrp (I \otimes \aut)$ are in $\mathbb{C}^3 \otimes \jotu$, namely they are matrix valued vectors with symmetric matrix elements;
\item the {\it sharp} product $\#$ of $\jotu$ matrices appearing in $\vspm \times \vrpm$ is a fundamental product in the theory of Jordan Algebras, \cite{McCrimmon}. It is the linearization of $\xs := x^2 - t(x) x - \frac12(t(x^2) - t(x)^2)I$, in terms of which we may write the fundamental cubic identity for $\jotn, n= 1,2,4,8$:
\be \xs\jdot  x = \frac13 t(\xs\!, x) I \quad \text{or} \quad x^3 - t(x) x^2 + t(\xs) x - \frac13 t(\xs\! , x) I = 0 \label{cubic} \ee
where $x^3 = x^2 \jdot x$ (notice that for $\joto$, because of non-associativity, $x^2 x \ne x x^2$ in general).
\end{enumerate}

The validity of the Jacobi identity for the algebra of matrices \eqref{mfq} with Lie product given by \eqref{fqcom} - \eqref{not1} derives from the Jacobi identity for $\rho(\fq)$ proven in \cite{T-2} together with Proposition \ref{sec:d(s)}.1, applied to $\ctm$ by trivially extending the three grading argument. The validity of the Jacobi identity, together with the fact that the representation $\rep$ fulfills the root diagram of $\ctm$ (as can be easily seen) proves that $\rep$ is indeed a representation of $\ctm$.

Before passing to $\esm$, let us point out that the cases of $\afm$ ($\mathbf{n}=2$) and $\dsm$ ($\mathbf{n}=4$) can be worked out in the same fashion as for $\ctm$, starting from the representations of $\es$ and $\est$ introduced in \cite{T-2}.


\section{$\mathbf{n}=8$ : Matrix representation of $\esm$}\label{Blues2}

We recall a few concepts and notations from \cite{T-2}. We use the notation $L_x z := x\jdot z$ and, for $\vx \in \mathbb{C}^3 \otimes \joto$ with components $(x_1, x_2, x_3)$, $L_\vx \in \mathbb{C}^3 \otimes L_{\joto}$ denotes the corresponding operator valued vector with components $(L_{x_1}, L_{x_2}, L_{x_3})$. We can write an element $\au$ of $\es$ as $\au =  L_x + \sum [L_{x_i},L_{y_i}]$ where $x,x_i,y_i \in \joto$ and $t(x) = 0$. The adjoint is defined by $\aud:= L_x - [L_{x_1},L_{x_2}]$. Notice that the operators $F := [L_{x_i},L_{y_i}]$ span the $\fq$ subalgebra of $\es$, the derivation algebra of $\joto$ . (Recall that the Lie algebra of the structure group of $\joto$ is $\es \oplus \mathbb{C}$.)\\

We remark that $(\au,-\aud)$ is a derivation in the Jordan Pair $(\joto,\jobto)$, and it is useful to recall that the relationship between the structure group of a Jordan algebra $J$ and the automorphism group of a Jordan Pair $V = (J,J)$ goes as follows, \cite{loos1}: if $g \in Str(J)$ then $(g, U^{-1}_{g(I)} g) \in Aut(V)$. In our case, for $g = 1 + \epsilon (L_x + F)$, at first order in $\epsilon$ (namely, in the tangent space of the corresponding group manifold) we get $ U^{-1}_{g(I)} g = 1 + \epsilon (- L_x + F) +O(\epsilon^2)$.\\

Next. we introduce a product $\mep$ such that $L_x \mep L_y :=  L_{x\cdot y} + [L_x, L_y]$, $F\mep L_x := 2 F L_x$ and $L_x \mep F :=2  L_x F$ for each component $x$ of $\vx \in \mathbb{C}^3\otimes \joto$ and $y$ of $\vy \in \mathbb{C}^3\otimes \joto$. If we denote by $[ ; ]$ the commutator with respect to the $\mep$ product, we also require that $[F_1 ; F_2] := 2 [F_1,F_2]$. We have that, $L_x \mep L_y + L_y \mep L_x = 2 L_{x\cdot y}$ and $[F; L_x] := F\mep L_x - L_x \mep F= 2 [F, L_x] = 2 L_{F(x)}$, where the last equality holds because $F$ is a derivation in $\joto$.\\

Therefore, for $\fqe \in \esm$, we write:
\begin{equation}
\rep(\fqe) = \left ( \begin{array}{cc} a\otimes Id + I \otimes \au & \Lvsp \\ \Lvsm & -I\otimes \aud
\end{array}\right)
\label{meo}
\end{equation}
where $a,\ \vspm$ are the same as in (\ref{sm}), $\au \in \es$, $I$ is the $3\times 3$ identity matrix, $Id := L_I$ is the identity operator in $L_{\joto}$: $L_I L_x= L_x$. Notice that $Id$ is the identity also with respect to the $\mep$ product.

By extending the $\mep$ product in an obvious way to the matrix elements \eqref{meo}, one achieves that $(I \otimes \au) \mep \Lvrp + \Lvrp \mep (I \otimes \aud) = 2 L_{(I \otimes \au) \vrp}$ and  $(I \otimes \aud) \mep \Lvrm + \Lvrm \mep (I \otimes \au) = 2 L_{(I \otimes \aud) \vrm}$.

After some algebra, the commutator of two matrices like \eqref{meo} can be computed to read :
\be
\begin{array}{c}
\left[
\left ( \begin{array}{cc} a\otimes Id + I \otimes \au & \Lvsp \\ \Lvsm & -I\otimes \aud
\end{array}\right) ,
\left ( \begin{array}{cc} b\otimes Id + I \otimes \bu &\Lvrp \\ \Lvrm & -I\otimes \bud
\end{array}\right) \right] \\  \\ :=
\left (\begin{array}{cc} C_{11} & C_{12}\\
C_{21} & C_{22}
 \end{array} \right), \hfill
\label{eocom}
\end{array}
\ee
where:
\be
\begin{array}{ll}
C_{11} &= [a\centerdot b] \otimes Id + 2 I \otimes [\au,\bu] + \Lvsp \diamond \Lvrm - \Lvrp \diamond \Lvsm \\ \\
C_{12} &=  (a \otimes Id) \Lvrp -  (b \otimes Id) \Lvsp +2  L_{(I \otimes \au) \vrp}\\
 &\phantom{:=} - 2 L_{(I \otimes \bu) \vsp} +  \Lvsm \times \Lvrm \\ \\
C_{21} &= - \Lvrm (a \otimes Id)  +  \Lvsm (b \otimes Id) - 2 L_{(I \otimes \aud) \vrm} \\
&\phantom{:=} +2  L_{(I \otimes \bud) \vsm} +  \Lvsp \times \Lvrp \\ \\
C_{22} &= 2 I \otimes [\aud,\bud] + \Lvsm \bullet \Lvrp - \Lvrm \bullet \Lvsp.
\end{array}
 \label{comreleo}
\ee
The products in \eqref{comreleo} are defined as follows :
\be
\begin{array}{ll}
\Lvsp \diamond \Lvrm &:= \left(\frac13 t(s^+_1, r^-_1) I - (1-(E_{23})_{ij}) t(s^+_i,r^-_j) E_{ij}\right) \otimes Id +\\
&\phantom{:=} I \otimes \left(\frac13 t(s^+_1, r^-_1) Id - L_{s^+_1 \cdot r^-_1} - [L_{s^+_1}, L_{r^-_1}] \right) \\ \\
\Lvsm \bullet \Lvrp &:= I \otimes (\frac13 t(s^-_1,r^+_1) Id - L_{s^-_1 \cdot r^+_1} - [L_{s^-_1}, L_{r^+_1}]) \\  \\
\Lvspm \times \Lvrpm &:= L_{\vspm \times \vrpm} = L_{\epsilon_{ijk} (s_j^\pm \# r_k^\pm)}
\end{array}
 \label{not1eo}
\ee

From the properties of the triple product of Jordan algebras, it holds that $ L_{s^+_1 \cdot r^-_1} + [L_{s^+_1}, L_{r^-_1}] = \frac12 V_{s^+_1 , r^-_1} \in \es\oplus \mathbb{C}$, \cite{T-2}. Moreover one can readily check that $[a_1^\dagger,b_1^\dagger] = - [a_1,b_1]^\dagger$ and $\Lvrm \bullet \Lvsp = I \otimes (\dfrac13 t(s^+_1,r^-_1) Id - L_{s^+_1 \cdot r^-_1} - [L_{s^+_1}, L_{r^-_1}])^\dagger$; this result implies that we are actually considering an algebra.

The validity of the Jacobi identity for the algebra of matrices \eqref{meo} with Lie product given by \eqref{eocom} - \eqref{not1eo} derives from the Jacobi identity\footnote{We would like to recall that the proof of the Jacobi identity given in \cite{T-2} strongly relies on identities deriving from the Jordan Pair axioms \cite{loos1}.} for $\rho(\eo)$ (proven in \cite{T-2}), together with Proposition \ref{sec:d(s)}.1, applied to $\esm$ by trivially extending the three grading argument. That the Lie algebra so represented is $\esm$ is made obvious by a comparison with the root diagram in figure \ref{fig:rootesm}.

\appendix

\section{Real Forms}

We use the notations of \cite{T-2}.
From the treatment in \cite{T-2}, a real form of octonions  is obtained by taking $\alpha_0^\pm, \alpha_k^\pm \in \rr$. The quaternionic subalgebra generated by $\rpm , \varepsilon_1^\pm$ is obviously a split form with nilpotent $\varepsilon_1^\pm$.\\

Another real form is obtained by taking complex coefficients with complex conjugation denoted by `$*$' subject to the conditions:
\be
\alpha_0^- = (\alpha_0^+)^* \ , \quad \alpha_1^- = -(\alpha_1^+)^* \ , \quad \alpha_3^- = (\alpha_2^+)^*
\ee
Its quaternionic subalgebra, generated by $1, u_7, iu_1, iu_4$, is also split with nilpotent $u_7 + i u_k \ , k=1,4$. It is equivalent to the one obtained with all real coefficients, which is generated by  $1, u_1, iu_4, iu_7$ upon cyclic permutation of the indices $7,4,1$.\\

Let us now restrict the Zorn matrix product \cite{T-2} to the sextonions and introduce the vectors
$$E_1 = (1,0,0) \qquad E_2 = (0,1,0) \qquad E_3 = (0,0,1)$$
  We get:
\begin{equation}
\begin{array}{c}
\left [ \begin{array}{cc} \alpha_0^+ & A^+ \\ A^- & \alpha_0^-
\end{array}\right]
\left [ \begin{array}{cc} \beta_0^+ & B^+ \\ B^- & \beta_0^-
\end{array}\right] \\ \\ =
\left [\begin{array}{cc} \alpha_0^+ \beta_0^+ - \alpha_1^+ \beta_1^- &
(\alpha_0^+ \beta_1^+ + \beta_0^- \alpha_1^+) E_1^+ \\
(\alpha_0^- \beta_1^- + \beta_0^+ \alpha_1^-)E_1^- & \alpha_0^- \beta_0^- -
\alpha_1^- \cdot \beta_1^+ \end{array} \right]\\ \\+
\left [\begin{array}{cc} 0 &
(\alpha_0^+ \beta_2^+ + \beta_0^- \alpha_2^+  + \alpha_3^- \beta_1^- -  \alpha_1^- \beta_3^-) E_2^+ \\
(\alpha_0^- \beta_3^- + \beta_0^+ \alpha_3^- + \alpha_1^+ \beta_2^+ -  \alpha_2^+ \beta_1^+)E_3^- & 0\end{array} \right]
\label{zorns}
\end{array}
\end{equation}

The algebra generated by $\rpm , \varepsilon_1^\pm$ is the quaternioc subalgebra. Its divisible real form is obtained by setting: $\alpha_0^- = (\alpha_0^+)^*$ and $\alpha_1^- = (\alpha_1^+)^*$ and it is a real linear span of $1, u_7, u_4, u_1$.\\
We now show that it is impossible to have e sextonion real algebra that has divisible quaternions as a subalgebra. To this aim we suppose  $\alpha_0^- = (\alpha_0^+)^*$ and $\alpha_1^- = (\alpha_1^+)^*$ and take in (\ref{zorns}) $\alpha_0^+ = \beta_0^+ = \beta_1^+ = 0$ - which implies $\alpha_0^- = \beta_0^- = \beta_1^- = 0$. The product (\ref{zorns}) shows that  the coefficients of $\varepsilon_2^+$ and $\varepsilon_3^-$ must be complex, hence each coefficient, say $\alpha_2^+$ contains 2 real parameters $a$ and $b$, and, in order to have a six dimensional real algebra, $\alpha_3^-$ viewed in $\rr_2$ must be a linear transformation $T$ of $(a,b)$, linearity being enforced by the linearity of the algebra. \\
We loosely write $\alpha_2^+ = T \alpha_3^-$. It is easy to show that $T^2 = Id$, namely $T$ is an involution. By playing with the coefficients in (\ref{zorns}) we can easily obtain $\alpha_2^+ = - T^2 \alpha_2^+ = - \alpha_2^+$  and similarly for $\alpha_3^-$, a contraddiction unless $\alpha_2^+ = \alpha_3^- = 0$.\\
This ends our proof.

\section*{Acknowledgements}

We would like to thank Leron Borsten and Bruce Westbury for useful
correspondence.


\end{document}